\documentclass[11pt]{amsart}
\usepackage{amsfonts,amssymb}
\usepackage{eucal}
\usepackage{graphics}
\usepackage{color}
\usepackage{amscd}

\usepackage[all,knot]{xy}
\usepackage{pstricks,}
\usepackage{pst-grad}
\usepackage{pst-plot}
\usepackage[tiling]{pst-fill}
\usepackage{graphics}

\newtheorem{Theorem}{Theorem}

\newtheorem{Conjecture}[Theorem]{Conjecture}

\newcommand{\C}{\mathbb C}
\newcommand{\R}{\mathbb R}
\renewcommand{\P}{\mathbb P}

\newcommand{\Z}{\mathbb Z}
\renewcommand{\S}{\mathbb S}
\newcommand{\Pic}{\textrm{Pic}}
\newcommand{\Sym}{\textrm{Sym}}
\newcommand{\Hilb}{\textrm{Hilb}}
\newcommand{\Diff}{\textrm{Diff}}
\newcommand{\Conf}{\textrm{Conf}}
\newcommand{\id}{\textrm{id}}

\newcommand{\rk}{\textrm{rk}}
\newcommand{\Kh}{\textrm{Kh}}

\newcommand{\Q}{\mathbb{Q}}

\newcommand{\Slice}{\EuScript{S}}

\newcommand{\Symp}{\textrm{Symp}}

\newcommand{\minushorizresolution}{\;\;\includegraphics{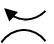}\;}
\newcommand{\plushorizresolution}{\;\;\includegraphics{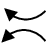}\;}
\newcommand{\minusvertresolution}{\;\;\includegraphics{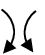}\;}
\newcommand{\plusvertresolution}{\;\;\includegraphics{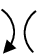}\;}
\newcommand{\minuscrossing}{\;\;\includegraphics{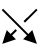}\;}
\newcommand{\pluscrossing}{\;\;\includegraphics{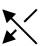}\;}

\title[Clay Institute Lectures]{Monodromy, vanishing cycles, knots \\ and the adjoint quotient}
\author{Ivan Smith}
\date{December 2004. Thanks to Aaron Lauda for help with the pictures. This work was partially supported by a grant NUF-NAL/00876/G from the Nuffield Foundation.}

\begin{document}
\maketitle

\begin{small}
These notes, adapted from two talks given at the Clay Institute Summer School on \emph{Floer homology, gauge theory, and low-dimensional topology} (Budapest, 2004), retain the format of the original lectures. The first part describes as background some of the geometry of symplectic fibre bundles and their monodromy.  The second part applies these general ideas to certain Stein fibre bundles that arise naturally in Lie theory, to construct an invariant of oriented links in the three-sphere (Section 2(h)).  Despite its very different origins, this invariant is conjecturally equal to the combinatorial homology theory defined by Mikhail Khovanov (Section 2(i)). In the hope of emphasising the key ideas, concision has taken preference over precision; there are no proofs, and sharp(er) forms of statements are deferred to the literature. 

Much of the first part I learned from, and the second part represents joint work with, Paul Seidel, whose influence and insights generously pervade all that follows.
\end{small}

\section{Monodromy, vanishing cycles}

Most of the material in this section is well-known; general references are \cite{McDS},\cite{Sei:LES},\cite{AurSmi}.  
\medskip

\noindent \textbf{1(A) Symplectic fibre bundles:} We will be concerned with fibrations $p:X\rightarrow B$ with symplectic base and fibre, or more precisely where $X$ carries a closed vertically non-degenerate 2-form $\Omega$, for which $d\Omega(u,v,\cdot)=0$ whenever $u$, $v$ are vertical tangent vectors.  If the fibration is proper, the cohomology class $[\Omega|_{\mathrm{Fibre}}]$ is locally constant, and parallel transport maps are symplectomorphisms.  Examples abound:

(1) A surface bundle over any space $\Sigma_g \rightarrow X \rightarrow B$ with fibre essential in homology can be given this structure; define $\Omega$ by picking any 2-form dual to the fibre and whose restriction to each fibre is an area form.  The homology constraint is automatically satisfied whenever $g\geq 2$ (evaluate the first Chern class of the vertical tangent bundle on a fibre).

(2) Given a holomorphic map $p:X \rightarrow B$ defined on a quasiprojective variety and which is smooth over $B^0 \subset B$, the restriction $p^{-1}(B^0)\rightarrow B^0$ defines a symplectic fibre bundle, where the 2-form $\Omega$ is the restriction of a K\"ahler form on $X$.  Such examples show the importance of \emph{singular fibres}.  Rational maps and linear systems in algebraic geometry provide a plethora of interesting singular fibrations.

(3) Contrastingly, the (singular) fibrations arising from moment maps, cotangent bundle projections and many dynamical systems have Lagrangian fibres, and fall outside the scope of the machinery we'll discuss.

Strictly, it is sensible to make a distinction between Hamiltonian and more general symplectic fibrations; essentially this amounts to the question of whether the vertically non-degenerate 2-form $\Omega$ has a closed extension to the total space, as in the cases above.  The subtlety will not play any role in what follows, but for discussion see \cite{McDS}.
\medskip

\noindent \textbf{1(B) Parallel transport:} A symplectic fibre bundle has a distinguished connexion, where the horizontal subspace at $x\in X$ is the symplectic orthogonal complement to the vertical distribution $\textrm{Hor}_x = \ker(dp_x)^{\perp_{\Omega}}$.  For K\"ahler fibrations, since the fibres are complex submanifolds, we can also define this as the orthogonal complement to $\ker(dp_x)$ with respect to the K\"ahler metric.  We should emphasise at once that, in contrast to the Darboux theorem which prevents local curvature-type invariants entering symplectic topology naively, there is no ``universal triviality'' result for symplectic fibre bundles.  The canonical connexion can, and often does, have curvature, and that curvature plays an essential part in the derivations of some of the theorems of the sequel.

Given a path $\gamma:[0,1]\rightarrow B$ we can lift the tangent vector $d\gamma/dt$ to a horizontal vector field on $p^{-1}(\gamma)$, and flowing along the integral curves of this vector field defines local \emph{symplectomorphisms} $h^{\gamma}$ of the fibres.

(1) If $p:X\rightarrow B$ is proper, the horizontal lifts can be globally integrated and we see $p$ is a fibre bundle with structure group $\Symp(p^{-1}(b))$.  Note that, since the connexion isn't flat, the structure group does not in general reduce to the symplectic mapping class group (of components of $\Symp$).

(2) Often there is a group $G$ acting fibrewise and preserving all the structure, in which case parallel transport will be $G$-equivariant.  An example will be given shortly.
\medskip

\noindent \textbf{1(C) Non-compactness:}
If the fibres are not compact, the local parallel transport maps may not be globally defined, since the solutions to the differential equations defining the integral curves may not exist for all times.  To overcome this, there are several possible strategies. The simplest involves estimating the parallel transport vector fields explicitly (which in turn might rely on choosing the right $\Omega$).  

Suppose for instance $p:X\rightarrow\C$ where $X$ has a K\"ahler metric, then for $V\in T_{p(x)}(\C)$ the lift $V^{hor} = V.\frac{(\nabla p)_x}{|\nabla p|^2}$.  If $p:(\C^n,\omega_{st}) \rightarrow \C$ is a homogeneous polynomial, clearly its only critical value is the origin, giving a fibre bundle over $\C^*$.  The identity $dp_x(x)=deg(p).p(x)$, together with the previous formula for $V^{hor}$, shows that the horizontal lift of a tangent vector $V \in T_p(x)(\C)=\C$ has norm $|V^{hor}| \leq \frac{|V|.|x|}{deg(p).|p(x)|}$.  On a fixed fibre $p=const$ this grows linearly with $|x|$ and can be globally integrated.

Corollary: For homogeneous polynomials $p:(\C^n,\omega_{st}) \rightarrow \C$ parallel transport is globally defined over $\C^*$.

Example: the above applies to the determinant mapping (indeed any single component of the characteristic polynomial or \emph{adjoint quotient} map), $\det: \textrm{Mat}_n(\C)\rightarrow\C$.  In this case, parallel transport is invariant under $SU_n\times SU_n$.  The monodromy of the associated bundle seems never to have been investigated.

For mappings $p:\C^n \rightarrow\C^m$ in which each component is a homogeneous polynomial but the homogeneous degrees differ, the above arguments do not quite apply but another approach can be useful. Since the smooth fibres are Stein manifolds of finite type, we can find a vector field $Z$ which points inwards on all the infinite cones.  By flowing with respect to a vector field $V^{hor}-\delta Z$, for large enough $\delta$, and then using the Liouville flows, we can define ``rescaled'' parallel transport maps $h_{resc}^{\gamma}: p^{-1}(t)\cap B(R) \hookrightarrow p^{-1}(t')$ on arbitrarily large pieces of a fixed fibre $p^{-1}(t)$, which embed such compacta symplectically into another fibre $p^{-1}(t')$.  This is not quite the same as saying that the fibres are globally symplectomorphic, but is enough to transport closed Lagrangian submanifolds around (uniquely to isotopy), and often suffices in applications.  For a detailed discussion, see \cite{SeiSmi}.  (In fact, if the Stein fibres are finite type and complete one can ``uncompress'' the flows above to show the fibres really are globally symplectomorphic, cf. \cite{KhoSei}.)
\medskip

\noindent \textbf{1(D) Vanishing cycles:} The local geometry near a singularity (critical fibre of $p$) shows up in the \emph{monodromy} of the smooth fibre bundle over $B^0$, i.e. the representation
$$\pi_1(B^0,b)\rightarrow \pi_0 (\Symp(p^{-1}(b),\Omega)).$$
Consider the ordinary double point (Morse singularity, node...) 
$$p:(z_1,\ldots,z_n)\mapsto \sum z_i^2.$$
The smooth fibres $p^{-1}(t)$, when equipped with the restriction of the flat K\"ahler form $(i/2)\sum_j dz_j\wedge d\overline{z}_j$ from $\C^n$, are symplectically isomorphic to $(T^*\S^{n-1}, \omega_{can})$.  Indeed, an explicit symplectomorphism can be given in co-ordinates by viewing 
$$T^*\S^{n-1} = \{ (a,b)\in\R^{n}\times\R^{n}\, | \, |a|=1, \langle a,b\rangle = 0\}$$
and taking $p^{-1}(1)\ni z \mapsto (\Re(z)/|\Re(z)|, -|\Re(z)|\Im(z))$. There is a distinguished Lagrangian submanifold of the fibre, the zero-section, which can also be defined as the locus of points which flow into the singularity under parallel transport along a radial line in $\C$. Accordingly, this locus -- which is $\{ (z_1,\ldots,z_n)\in\R^n \, | \, \sum_i |z_i|^2 = 1 \}$ in co-ordinates -- is also called the \emph{vanishing cycle} of the singularity.

Lemma: The monodromy about a loop encircling $0 \in \C$ is a Dehn twist in the vanishing cycle.

To define the Dehn twist, fix the usual metric on $\S^{n-1}\subset \R^n$, which identifies $T^*\S^{n-1} \cong T\S^{n-1}$. The Dehn twist is the composite of the time $\pi$ map of the geodesic flow on the unit disc tangent bundle $U(T\S^{n-1})$ with the map induced by the antipodal map; it's antipodal on the zero-section and vanishes on the boundary $\partial U$.  If $n=2$ this construction is classical, and we get the usual Dehn twist on a curve in an annulus $T^*\S^1$.
\medskip

\noindent \textbf{1(E) Variant constructions:} There are two useful variants of the above model:

(1a) The relative version of the above geometry: If we have $p:X\rightarrow \C$ and the fibre $X_0$ over $0\in\C$ has smooth singular locus $Z^c$, with normal data locally holomorphically modelled on the map $\sum_{i=1}^n z_i^2$, we say $p$ has a ``fibred $A_1$-singularity''.  Then the nearby smooth fibre $X_t$ contains a relative vanishing cycle $\S^{n-1}\rightarrow Z \rightarrow Z^c$. An open neighbourhood $U$ of $Z\subset X_t$ is of the form $T^*\S^{n-1}\hookrightarrow U \rightarrow Z^c$, and the monodromy about $0$ is a ``fibred Dehn twist'' (the above construction in every $T^*$-fibre).

(1b) Examples: Fix a stable curve over a disc $f:\mathcal{X}\rightarrow\Delta$ with all fibres smooth except for an irreducible curve with a single node over the origin.  The relative Picard fibration $\Pic(f)\rightarrow\Delta$ has a singular fibre over $0$ with a fibred $A_1$-singularity, and with singular set the Picard of the normalisation $\tilde{C}$ of $C=f^{-1}(0)$.  The relative Hilbert scheme $\Hilb^r(f)\rightarrow\Delta$ has a fibred $A_1$-singularity over $0$ with singular locus $\Hilb^{r-1}(\tilde{C})=\Sym^{r-1}(\tilde{C})$.  In both these cases $n=2$.  For the relative moduli space of stable rank two bundles with fixed odd determinant, there is again a model for the compactification (symplectically, not yet constructed algebraically) with a fibred $A_1$-singularity, but this time with $n=4$ and hence $\S^3\cong SU(2)$ vanishing cycles \cite{Sei:notes}.

(2a) Morsification: If we have a ``worse'' isolated singular point at the origin of a hypersurface defined e.g. by some polynomial $P(x)=0$ then often we can perturb to some $P_{\varepsilon}(x)=0$ which is an isotopic hypersurface outside a compact neighbourhood of the original singularity but which now only has a collection of finitely many nodes.  In particular, the global monodromy of the projection of $P^{-1}(0)$ to the first co-ordinate in $\C$ about a large circle is a \emph{product of Dehn twists in Lagrangian spheres} in the generic fibre.

(2b) Examples: The triple point $x^3+y^3+z^3+t^3=0$ can be perturbed to a hypersurface with 16 nodes, or completely smoothed to give a configuration of 16 Lagrangian vanishing cycles, cf. \cite{SmiTho}.

The second result enters into various ``surgery theoretic'' arguments, along the following lines. Given a symplectic manifold with a tree-like configuration of Lagrangian spheres which matches the configuration of vanishing cycles of the Morsified singularity, we can cut out a (convex) neighbourhood of the tree and replace it with the resolution (full blow-up) of the original singular point, and this is a \emph{symplectic} surgery.  Examples are given in \cite{SmiTho}. 
\medskip

\noindent \textbf{1(F) Lefschetz fibrations:}
A remarkable theorem due to 
Donaldson asserts that every symplectic manifold admits a Lefschetz pencil.  In dimension four, this comprises a map $f:X\backslash\{b_i\} \rightarrow \S^2$ submersive away from a finite set $\{p_j\}$, and with $f$ given by $z_1/z_2$ near $b_i$ and $z_1z_2$ near $p_j$.  Removing fibres $f^{-1}(p_j)$ gives a symplectic fibre bundle over $\S^2\backslash\{f(p_j)\}$, and the global monodromy is encoded as a word in positive Dehn twists in $\Gamma_g = \pi_0(\Symp(\Sigma_g))$.  In general, Donaldson's theory of symplectic linear systems reduces a swathe of symplectic topology to combinatorial group theory, and places issues of monodromy at the centre of the symplectic stage.

Example: The equation 
$$\left( \left( \begin{array}{cc} 1&1 \\ 0&1\end{array}\right)\left(\begin{array}{cc} 1&0\\-1&1\end{array}\right) \right)^{6n} \ = \ (AB)^{6n} \ =\ I$$ 
encodes the elliptic surface $E(n)$ as a word in $SL(2,\Z)$.  The fact that all words in matrices conjugate to $A\in SL_2(\Z)$ are equivalent by the Hurwitz action of the braid group to a word of this form, for some $n$, gives an algebraic proof that the $E(n)$ exhaust all elliptic Lefschetz fibrations, hence all such are K\"ahler.  In turn, from this one can deduce that degree 4 symplectic surfaces in $\C\P^2$ are isotopic to complex curves.

(1) This kind of algebraic monodromy encoding generalises branched covers of Riemann surfaces and gives (in principle) a classification of integral symplectic 4-manifolds.

(2) The importance of \emph{Lagrangian} intersection theory -- i.e. geometric and not algebraic intersections of curves on a 2d surface -- already becomes clear. 

(3) Donaldson has suggested that the algebraic complexity of Lefschetz fibrations might be successfully married with the algebraic structure of Floer homology \cite{Do:PVCFH}.  Steps in this direction were first taken by Seidel in the remarkable \cite{Sei:VM}, see also \cite{Sei:K3},\cite{AKO}.

Donaldson's initial ideas have been developed and extended in a host of useful and indicative directions: we mention a few.  Lefschetz pencils can be constructed adapted to embedded symplectic submanifolds or Lagrangian submanifolds \cite{AMP} (in the latter case one extends a Morse function on $L$ to a Lefschetz pencil on $X\supset L$); there are higher-dimensional linear systems, leading to iterative algebraic encodings of symplectic manifolds \cite{Aur:TJ}; analogues exist in contact topology \cite{Pre} and, most recently, for (non-symplectic) self-dual harmonic 2-forms on four-manifolds \cite{ADK}.  In each case, the techniques give an algebraic encoding of some important piece of geometric data.

Challenge: show that symplectomorphism of integral simply-connected symplectic 4-manifolds is (un)decidable.
\medskip

\noindent \textbf{1(G) Counting sections:} A good way to define an invariant for a Lefschetz fibration is to replace the fibres with something more interesting and then count holomorphic sections of the new beast.  In other words, one studies the Gromov-Witten invariants for those homology classes which have intersection number 1 with the fibre of the new fibration.  

Explicitly, suppose we have a moduli problem on Riemann surfaces in the following sense: $\Sigma \mapsto \mathcal{M}(\Sigma)$ associates to a Riemann surface $\Sigma$ some projective or quasiprojective moduli space, with a relative version for families of irreducible stable curves.  A Lefschetz fibration $f:X\rightarrow\S^2$ gives rise to a relative moduli space $F:\mathcal{M}(f)\rightarrow\S^2$; our assumptions on the moduli problem should ensure that this is smooth and symplectic, and is either convex at infinity, compact or has a natural compactification.  Then we associate $(X,f)\mapsto Gr_A(F)$ where $A$ is the homology class of some fixed section.  This follows a philosophy derived from algebraic geometry: holomorphic sections of a family of moduli spaces on the fibres should be ``equivalent information'' to data about geometric objects on the total space which could, in principle, be defined without recourse to any given fibration structure.  Naive as this sounds, the theory is not entirely hopeless in the actual examples.

(1a) For $f:X^4 \rightarrow \S^2$, replace $X_t$ by $\Sym^r(X_t)$ and desingularise, forming the relative Hilbert scheme \cite{DonSmi} to get $F:X_r(f) \rightarrow \S^2$. Obviously sections of this new fibration are related to 2-cycles in the original four-manifold.  A pretty theorem due to Michael Usher \cite{Ush} makes this intuition concrete and sets the theory in a very satisfactory form:  the Gromov-Witten invariants $\mathcal{I}_{X,f}$ counting sections of $F$, known as the \emph{standard surface count}, are equal to Taubes' $Gr(X)$.  In particular, the invariants are independent of $f$, as algebraic geometers would expect.  

(1b) Application \cite{DonSmi}, \cite{Smi:ST}: if $b_+>2$ the invariant $\mathcal{I}_{X,f}(\kappa) = \pm 1$, where $\kappa$ refers to the unique homology class of section for which the cycles defined in $X$ lie in the class Poincar\'e dual to $K_X$. This gives a Seiberg-Witten free proof of the fact that, for minimal such manifolds, $c_1^2(X)\geq 0$.  The key to the argument is the Abel-Jacobi map $\Sym^r{\Sigma}\rightarrow\Pic^r(\Sigma)\cong\mathbb{T}^{2g}$, which describes $\Sym^r(\Sigma)$ as a family of projective spaces over a torus; for the corresponding fibrations with fibre $\P^n$ or $\mathbb{T}^{2g}$ one can compute moduli spaces of holomorphic sections explicitly, and hence compute Gromov invariants.

(2) One can also count sections of symplectic Lefschetz fibrations over surfaces with boundary, provided suitable Lagrangian boundary conditions are specified.  In place of absolute invariants one obtains invariants living in Floer homology groups associated to the boundary, or formulated differently morphisms on Floer homology groups. This is reminiscent of the formalism of Topological Quantum Field Theory;  such ideas are central to the main theorem of \cite{Sei:LES}.

In (1), the fact that the compactifications of the relative moduli spaces exist and are smooth can be understood in terms of the local geometry of fibred $A_1$-singularities and normal crossings, as in the discussion of Section 1(E), (1b) above. 
\medskip

\noindent \textbf{1(H) Braid relations:}
It's harder to get invariants of the total space straight out of the monodromy of a fibre bundle, but it is very natural to study $\pi_0\Symp(Fibre)$ this way.
Let $p:\C^{n+1}\rightarrow\C$ be given by $\{x^{k+1}+\sum_{j=1}^n y_j^2 = \varepsilon$\}. There are $(k+1)$ critical values, and if we fix a path between two then we can construct a Lagrangian $\S^{n+1}$ in the total space by ``matching'' vanishing $\S^n$-cycles associated to two critical points \cite{Sei:K3}.  This is just the reverse process of finding a Lefschetz fibration adapted to a given Lagrangian $(n+1)$-sphere, by extending the obvious Morse function from $\S^{n+1}$ to the total space, mentioned above.

Lemma: For two Lagrangian spheres $L_1$, $L_2$ meeting transversely in a point, the Dehn twists $\tau_{L_i}$ satisfy the braid relation $\tau_{L_1}\tau_{L_2}\tau_{L_1}=\tau_{L_2}\tau_{L_1}\tau_{L_2}$.

The proof of this is by direct computation \cite{Sei:knot}. In the lowest dimensional case $n=1$ it is completely classical.  A disguised version of the same Lemma will underly central properties of a fibre bundle of importance in our application to knot theory in the second part.

Corollary: If $X$ contains an $A_k$-chain of Lagrangian spheres, there is a natural homomorphism $Br_k \rightarrow \pi_0(\Symp_{ct}(X))$.

These homomorphisms have come to prominence in part because of mirror symmetry, cf. \cite{SeiTho}.
The relevant chains of Lagrangian spheres can be obtained by Morsifying $A_{k+1}$-singularities.  The existence of the homomorphism, of course, gives no information on its non-triviality; we address that next.
\medskip

\noindent \textbf{1(I) Simultaneous resolution:}
The map $\C^3\rightarrow \C$ which defines a node $z_1^2+z_2^2+z_3^2$ has fibres $T^*\S^2$ and one singular fibre. If we pull back under a double cover $\C \rightarrow \C, w \mapsto w^2$ then we get $\sum z_i^2 -w^2=0$, i.e. a 3-fold ODP, which has a \emph{small resolution}; replacing the singular point by $\C P^1$ gives a smooth space.

Corollary: the fibre bundle upstairs is \emph{differentiably} trivial, since it completes to a fibre bundle over the disc.  

Seidel showed in \cite{Sei:DT} that this is not true symplectically; the Dehn twist in $T^*\S^2$ has infinite order as a symplectomorphism.  So the natural map 
$$\pi_0(\Symp_{ct}(T^*\S^2))\rightarrow\pi_0(\Diff_{ct}(T^*\S^2))$$ 
has infinite kernel -- the interesting structure is \emph{only visible symplectically}. For a smoothing of the $A_{k+1}$-singularity above, a similar picture shows the braid group acts \emph{faithfully} by symplectomorphisms but factors through $\Sym_k$ acting by diffeomorphisms (compactly supported in each case). The injectivity is established by delicate Floer homology computations \cite{KhoSei}.  Such a phenomenon is at least possible whenever one considers families with \emph{simultaneous resolutions}; that is, a family $\mathcal{X}\stackrel{\phi}{\longrightarrow} B$ for which there is a ramified covering $\tilde{B}\rightarrow B$ and a family $\tilde{\mathcal{X}}\stackrel{\tilde{\phi}}{\longrightarrow}\tilde{B}$ with a map $\pi:\tilde{\mathcal{X}}\rightarrow \mathcal{X}$ and with
$$\pi:\tilde{\mathcal{X}}_t = \tilde{\phi}^{-1}(t) \rightarrow \phi^{-1}(\pi(t)) = \mathcal{X}_t$$
a resolution of singularities for every $t \in \tilde{B}$.  The small resolution of the 3-fold node will be the first in a sequence of simultaneous resolutions considered in the second section, and in each case the inclusion of $\Symp_{ct}$ into $\Diff_{ct}$ of the generic fibre will have infinite kernel.
\medskip

\noindent \textbf{1(J) Long exact sequences:}
Aside from their role in monodromy, Lagrangian spheres and Dehn twists also give rise to special structures and properties of Floer cohomology. Suppose $L_1,L_2$ are Lagrangians in $X$ and $L\cong \S^n$ is a Lagrangian sphere.  The main theorem of \cite{Sei:LES} is the following:

Theorem: (Seidel) Under suitable technical conditions, there is a long exact triangle of Floer cohomology groups
\begin{small}
$$HF(L_1,L_2)\rightarrow HF(L_1,\tau_L(L_2))\rightarrow HF(L_1,L)\otimes HF(L,L_2).$$
\end{small}

The technical conditions are in particular valid for exact Lagrangian submanifolds of a Stein manifold of finite type; in this setting there is no bubbling, and the manifold will be convex at infinity which prevents loss of compactness from solutions escaping to infinity.  Hence,  the Floer homology groups are well-defined; if moreover the Stein manifold has $c_1=0$ (for instance is hyperk\"ahler), the groups in the exact triangle can be naturally $\Z$-graded. 

Corollary (\cite{Sei:FD}, Theorem 3): For a Lefschetz pencil of K3 surfaces in a Fano 3-fold, the vanishing cycles $\{L_j\}$ ``fill'' the generic affine fibre: every closed Lagrangian submanifold disjoint from the base locus and with well-defined Floer homology must hit one of the $\{L_j\}$.  

Proof: The global monodromy acts as a \emph{shift} on (graded) $HF^*$, so if $K$ is disjoint from all the spheres then the exact sequence shows 
$HF(K,K) = HF(K,K)[shifted]$.  Iterating, and recalling that $HF(K,K)$ is supported in finitely many degrees, this forces $HF(K,K)=0$. But this is impossible for any homologically injective Lagrangian submanifold, by general properties of Floer theory, cf. \cite{FO3},  which completes the contradiction.  

There are simpler proofs that any Lagrangian must intersect one of the vanishing cycles, but this gives a bit more: the vanishing cycles ``generate'' Donaldson's quantum category of the K3 (the underlying homological category of the Fukaya category). The Corollary above was in part motivated by an older and easier result, specific to the situation for curves in Riemann surfaces, given in \cite{Smi:GMHD}. 

In the second part we will focus attention on a Stein manifold $Y_m$ which also contains a distinguished finite collection of Lagrangian submanifolds (cf. Section 2(f) below), which conjecturally generate the quantum category of $Y_m$ in a similar way.  However, these arise not as vanishing cycles of a pencil but from the components of a ``complex Lagrangian'' small resolution, giving another point of contact between the two general themes of the last section.


\section{Knots, the adjoint quotient}

All the material of this section is joint work with Paul Seidel; we were considerably influenced by the ideas of Mikhail Khovanov.  References are \cite{SeiSmi},\cite{SeiSmi2},\cite{Kho}.

\noindent \textbf{2(a) Knot polynomials:} The Jones polynomial and Alexander polynomial $V_K(t), \Delta_K(t)$ are powerful knot invariants defined by \emph{skein relations}.  They are Laurent polynomials in $t^{\pm 1/2}$ determined by saying $V_K(U)=1=\Delta_K(U)$, for $U$ an unknot, and also that
$$t^{-1}V_{L_+} - tV_{L_-} + (t^{-1/2}-t^{1/2})V_{L_0} = 0;$$
$$\Delta_{L_+} - \Delta_{L_-} - (t^{-1/2}-t^{1/2})\Delta_{L_0} = 0.$$
Here the relevant links differ only near a single crossing, where they look as in the picture below:
\[
\xy
   (6,6)*{}="tl";
   (-6,6)*{}="tr";
   (6,-6)*{}="bl";
   (-6,-6)*{}="br";
        {\ar|{\hole \; \hole \; \hole \; \hole \; \hole \; \hole }  "bl";"tr"};
        {\ar  "br";"tl"};
   (0,-10)*{L_{+}};
 \endxy
\qquad \qquad
 \xy
   (-6,6)*{}="tl";
   (6,6)*{}="tr";
   (-6,-6)*{}="bl";
   (6,-6)*{}="br";
        {\ar|{\hole \; \hole \; \hole \; \hole \; \hole \; \hole }  "bl";"tr"};
        {\ar  "br";"tl"};
   (0,-10)*{L_{-}};
 \endxy
 \qquad \qquad
 \xy
   (-6,6)*{}="tl";
   (6,6)*{}="tr";
   (-6,-6)*{}="bl";
   (6,-6)*{}="br";
        "bl";"tl" **\crv{(-1,0)}?(1)*\dir{>};
        "br";"tr" **\crv{(1,0)}?(1)*\dir{>};
   (0,-10)*{L_{0}};
 \endxy
\]
The Alexander polynomial is well-understood geometrically, via homology of an infinite cyclic cover $H_1(\widetilde{\S^3\backslash K})$ \cite{Lic}.  There is also an interpretation in terms of 3-dimensional Seiberg-Witten invariants, beautifully explained in \cite{Do:MT}.  The Jones polynomial is more mysterious, although it does have certain representation theoretic incarnations in the theory of quantum groups and loop groups. The Jones polynomial solved a host of conjectures immediately after its introduction, one famous one being the following:

Example: (Kaufmann) A connected reduced alternating diagram for a knot exhibits the minimal number of crossings of any diagram for the knot.  [Reduced: no crossing can be removed by ``flipping'' half the diagram.]

Before moving on, it will be helpful to rephrase the skein property in the following slightly more involved fashion.
\begin{equation} \label{eq:skein}
\begin{aligned}
 &
 t^{-1/2} V_{\plushorizresolution}
 + t^{3v/2} V_{\plusvertresolution}
 + t^{-1} V_{\pluscrossing} = 0, \\
 &
 t^{3v/2} V_{\minushorizresolution}
 + t^{1/2} V_{\minusvertresolution}
 + t V_{\minuscrossing} = 0.
\end{aligned}
\end{equation}
Here $v$ denotes the signed number of crossings between the arc ending at the top left of the crossing and the other connected components of the diagram. Some of the arcs have no labelled arrow since resolving a crossing in one of the two possible ways involves a non-local change of orientation, but the relations are between the polynomials of \emph{oriented} links. Obviously these two equations together imply the original skein relation.  (The Jones polynomial of a \emph{knot} is independent of the choice of orientation, but for links this is no longer true.)
\medskip

\noindent \textbf{2(b) Khovanov homology:}
Mikhail Khovanov (circa 1998) ``categorified'' the Jones polynomial -- he defines combinatorially an invariant $K \mapsto Kh^{*,*}(K)$ which is a $\Z\times\Z$-graded abelian group, and such that

(i) $Kh^{0,*}(U_n) = H^*((\S^2)^n)[-n]$, where $U_n$ is an n-component unlink (and the cohomology is concentrated in degrees $(0,*)$);

(ii) Skein-type exact sequence: for oriented links as indicated, there are long exact sequences which play the role of 
\eqref{eq:skein} above:
\begin{multline}
 \label{eq:plustriangle}
 \qquad \cdots \longrightarrow \Kh^{i,j}(\!\pluscrossing) \longrightarrow
 \Kh^{i,j-1}(\!\plushorizresolution) \longrightarrow
 \Kh^{i-v,j-3v-2}(\plusvertresolution) \qquad\qquad \\ \longrightarrow
 \Kh^{i+1,j}(\!\pluscrossing) \longrightarrow \cdots \qquad \qquad \qquad
 \qquad \quad \qquad \qquad \qquad
\end{multline}
and
\begin{multline}
 \label{eq:minustriangle}
 \quad \cdots \longrightarrow  \Kh^{i,j}(\minuscrossing) \longrightarrow
 \Kh^{i-v+1,j-3v+2}(\minushorizresolution) \longrightarrow
 \Kh^{i+1,j+1}(\minusvertresolution) \quad\quad \\ \longrightarrow
  \Kh^{i+1,j}(\minuscrossing) \longrightarrow \cdots \qquad \qquad \qquad
  \qquad \quad \qquad \qquad \qquad
\end{multline}

(iii) As an easy consequence of (ii), a change of variables recovers Jones:
$$\frac{1}{q+q^{-1}}\, \sum_{i,j} (-1)^i \rk_{\Q}\Kh^{i,j}(K)q^j \ = \ V_K(t)|_{q=-\sqrt{t}}.$$

Note the exact sequences are not quite skein relations, since they do not figure the crossing change, rather than the two different crossing resolutions (sometimes called the horizontal and vertical resolutions, as in the next picture).
\[
 \xy
   (-6,6)*{}="tl";
   (6,6)*{}="tr";
   (-6,-6)*{}="bl";
   (6,-6)*{}="br";
        "tl";"tr" **\crv{(0,1)};
        "bl";"br" **\crv{(0,-1)};
   (0,-10)*{L_{{\rm hor}}};
 \endxy
 \qquad \qquad
 \xy
   (6,6)*{}="tl";
   (-6,6)*{}="tr";
   (6,-6)*{}="bl";
   (-6,-6)*{}="br";
        {\ar@{-}|{\hole \; \hole \; \hole \; \hole \; \hole \; \hole }  "bl";"tr"};
        {\ar@{-}  "br";"tl"};
   (0,-10)*{L_{\rm cross}};
 \endxy
 \qquad \qquad
  \xy
   (-6,6)*{}="tl";
   (6,6)*{}="tr";
   (-6,-6)*{}="bl";
   (6,-6)*{}="br";
        "bl";"tl" **\crv{(-1,0)};
        "br";"tr" **\crv{(1,0)};
   (0,-10)*{L_{\rm vert}};
 \endxy
\]
\medskip

Khovanov homology is known to be a strictly stronger invariant than the Jones polynomial, but its principal interest lies in its extension to a ``Topological Quantum Field Theory''; cobordisms of knots and links induce canonical homomorphisms of Khovanov homology.  Relying heavily on this structure, at least one beautiful topological application has now emerged:

Example: Rasmussen \cite{Ras} uses $Kh^{*,*}$ to compute the unknotting number, which is also the slice genus, of torus knots, $\textrm{Unknot}(T_{p,q}) = (p-1)(q-1)/2$.  

This result, first proved by Kronheimer and Mrowka, was formerly accessible only via adjunction-type formulae in gauge theory (or the rebirth of gauge theory via Ozsvath and Szabo); by all comparison, Rasmussen's combinatorial proof represents an enormous simplification.  One current limitation on Khovanov homology is precisely that its mystery makes it unclear which, comparable or other, problems it could profitably be applied to.
\medskip

\noindent \textbf{2(c) Invariants of braids:}  Here is a general way to (try to) define knot invariants using symplectic geometry. We begin with: 

(1) a symplectic fibre bundle $\mathcal{Y}\rightarrow\Conf_{2n}(\C)$ over the configuration space of unordered $2n$-tuples of points in $\C$.  Suppose parallel transport is well-defined, or at least its rescaled cousin from Section 1(B).

(2) a distinguished (to isotopy) Lagrangian submanifold $L\subset \mathcal{Y}_t$ in some distinguished fibre over $t\in \Conf_{2n}(\C)$.

Given a braid $\beta$ on $2n$-strands, i.e. a loop in the base, we can use parallel transport to get Lagrangian submanifolds $L, \beta(L)\subset \mathcal{Y}_t$ and then consider the Lagrangian Floer homology group $\beta \mapsto HF(L,\beta(L))$. This is the homology of a chain complex generated by intersection points, with boundary maps defined by counting pseudoholomorphic discs with boundary on the Lagrangian submanifolds as in the picture below.

\[
\begin{pspicture}(2,1.5)
 \pscircle[linecolor=black, fillcolor=lightgray,
 fillstyle=gradient,
  gradbegin=lightgray, gradend=darkgray, gradmidpoint=1,gradangle=330](0,0){1}
  \psdot(0,.99) \psdot(0,-.99)  \rput(-1.1,.5){$\ell^{+}$} \rput(1.2,.5){$\ell^{-}$}
  \rput(1.2,-.5){$D^2$}
 \end{pspicture}
\qquad \xy {\ar (-5,0);(5,0)}; \endxy \qquad \qquad \qquad \qquad
\psset{unit=0.5cm}
\begin{pspicture}(2,1.5)
\begin{psclip}{ \pscustom[linestyle=none,linecolor=white]{
  \psbezier(0,0)(1.5,-0.2)(2.5,1)(2,2.5)\lineto(-1,2.5)}
} \pscustom[linestyle=none, fillcolor=lightgray,
fillstyle=gradient,gradbegin=lightgray,
gradend=darkgray,gradmidpoint=1,gradangle=330]{
 \psbezier(0,0)(0,1)(1.5,2)(2.8,2)\lineto(2.8,-1)}
 \end{psclip}
 \psbezier(0,0)(0,1)(1.5,2)(2.8,2)
 \psbezier(0,0)(1.5,-0.2)(2.5,1)(2,2.5)
 \psbezier(0,0)(0,-1)(-1.5,-2)(-2.8,-2)
 \psbezier(0,0)(-1.5,0.2)(-2.5,-1)(-2,-2.5)
  \rput(2.4,1){$\ell^{+}$}  \rput(.9,1.7){$\ell^{-}$}
   \rput(-2,-2.9){$L_2$} \rput(-3.2,-2){$L_1$}
   \psdot(2.12,1.9) \psdot(0,0)
\end{pspicture}
\]
\bigskip
\bigskip

Caution: we're ignoring all technical difficulties.  As before, well-definition of Floer homology relies on overcoming compactness problems, but for exact Lagrangians in finite type Stein manifolds this is standard.  If the Lagrangians are spin, there are coherent orientations and Floer homology can be defined with $\Z$-coefficients.  If the Lagrangians have $b_1=0$ (so zero Maslov class) and the ambient space has trivial first Chern class, the Lagrangians can be graded and Floer homology will be $\Z$-graded.

In the discussion so far, we could obtain invariants of braids on any number of strands. The restriction to the even-strand case comes in making the connection to the theory of knots and links, which we do below.
\medskip

\noindent \textbf{2(d) Markov moves:}  It is well-known that every oriented link can be obtained as the ``closure'' of a braid in the fashion given in the following diagram: one goes from $Br_n \ni \beta \mapsto \beta\times\id \in Br_{2n}$ and then caps off top and bottom with a collection of nested horseshoes.  Such a representation of oriented links is enormously non-unique, but the equivalence relation on braids that generates this non-uniqueness is well-understood, and generated by the so-called Markov moves. The first is conjugation $\beta \mapsto \sigma\beta\sigma^{-1}$ by any $\sigma \in Br_n$, and the second -- which is more interesting, since it changes the number of strands of the braid -- involves linking in an additional strand by a single positive or negative half-twist, giving the $II^+$ and $II^-$ stabilisations. All are pictured below.  (To see that the link closure is canonically oriented, put an ``upwards'' arrow on all the parallel right hand strands.)
\bigskip

(1) Link closure:

\[
 \beta
 \qquad
  \xy
   (0,0)*+{\includegraphics{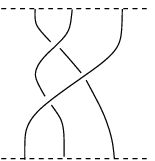}};
   (0,5)*{};
 \endxy
 \quad
 \mapsto
 \quad
 \xy
   (0,0)*+{\includegraphics{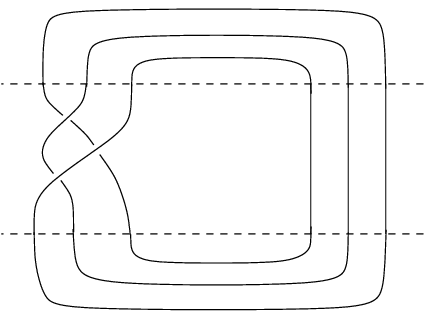}};
 \endxy
\]

(2) Markov $I$:

\[
 \beta
 \qquad
  \xy
   (0,0)*+{\includegraphics{R3-30.eps}};
 \endxy
 \quad
 \mapsto
 \quad
 \xy
   (0,0)*+{\includegraphics{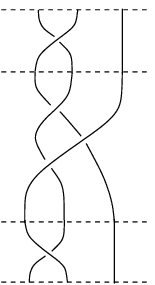}};
 \endxy
\]

(3) Markov $II^+$:

\[
 \beta
 \qquad
  \xy
   (0,0)*+{\includegraphics{R3-30.eps}};
 \endxy
 \quad
 \mapsto
 \quad
 \xy
   (0,0)*+{\includegraphics{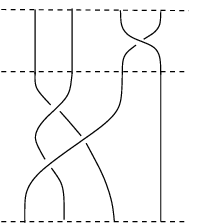}};
 \endxy
\]

It follows that if the association $\beta \mapsto HF(L,(\beta\times\id)(L))$ of Section 2(c) is invariant under the Markov moves, then it in fact defines an invariant of oriented links in the three-sphere.  We will now turn to a particular case in which precisely this occurs.
\medskip

\noindent \textbf{2(e) The adjoint quotient:} 
We will get our family of symplectic manifolds from (a cousin of) the characteristic polynomial mapping, also called the adjoint quotient $\chi: \frak{sl}_m \rightarrow \C^{m-1}$ which is smooth over $\Conf^0_{m}(\C)$, the space of balanced configurations, i.e. symmetric functions of distinct eigenvalues of trace-free matrices.  The following is the content of the Jacobson-Morozov theorem:

Fact: given a nilpotent matrix $N^+\in\frak{sl}_m(\C)$, there is a unique conjugacy class of homomorphisms $\rho:\frak{sl}_2\rightarrow\frak{sl}_m$ such that \begin{small} $\left(\begin{array}{cc} 0&1\\0&0\end{array}\right) \stackrel{\rho}{\mapsto}N^+$ \end{small}  

Let $N^-$ be the image of the other standard nilpotent $\left( \begin{array}{cc}0&0\\1&0\end{array}\right)$ in $\frak{sl}_2$. Then $N^+ + \ker(ad \, N^-)$ is a ``transverse slice'' to the adjoint action, i.e. it's an affine subspace of $\frak{sl}_m$ which intersects the adjoint orbit of $N^+$ inside the nilpotent cone only at $N^+$.  (Given the JM theorem, this is an easy fact about $\frak{sl}_2$-representations.)

Lemma: the restriction of $\chi$ to such a transverse slice $\Slice$ is still a fibre bundle over $\Conf^0_m(\C)$.

For a suitable K\"ahler form, the rescaled parallel transport construction for this fibre bundle can be pushed through, and closed Lagrangian submanifolds transported into any desired fibres. This follows the general programme outlined at the end of Section 1(B).  Although we will not dwell on the details here, we should say at once that the relevant symplectic forms are \emph{exact}, and are not related to the Kostant-Kirillov forms that also arise when dealing with the symplectic geometry of adjoint orbits.
\medskip

\noindent \textbf{2(f) Simultaneous resolutions:}
Grothendieck gave a simultaneous resolution of $\chi|\Slice$: replace a matrix $A$ by the space of pairs $(A,\mathcal{F})$ where $\mathcal{F}$ is a flag stabilised by $A$.  This orders the eigenvalues, i.e. the resolution involves base-changing by pulling back under the symmetric group.  Hence, via Section~1(I), the differentiable monodromy of the fibre bundle $\chi:\Slice \rightarrow \C^{m-1}$ factors through $\Sym_m$; we get a diagram as follows, writing $Y_t$ for a fibre of the slice over some point $t$:
\[
\begin{CD}
 Br_{m} = \pi_1(\Conf_{m}^0(\C)) \hspace{-1em} @>{ \qquad}>> \pi_0(\text{Symp}(Y_{t}))  \\
 @VVV @VVV \\
 \Sym_{m}  @>{\qquad\qquad\!\!\!}>> \pi_0(\text{Diff}(Y_{t}))
\end{CD}
\]
All representations of symmetric groups (and more generally Weyl groups) arise this way, in what is generally known as the \emph{Springer correspondence}. As in Section~1(I), the symplectic monodromy is far richer (perhaps even faithful?).

Example: $\frak{sl}_2(\C)$ and $N^+=0$ so $\Slice=\frak{sl}_2$; then $\chi$ is the map $(a,b,c,-a)\mapsto -a^2-bc$ which after a change of co-ordinates is the usual node, with generic fibre $T^*\S^2$.

Example: $\frak{sl}_{2m}(\C)$ and $N^+$ with two Jordan blocks of equal size. Then the slice $\Slice_m \cong \C^{4m-1}$ is all matrices of $2\times 2$-blocks with $I_2$ above the diagonal, any $(A_1,\ldots,A_m)$ in the first column with $tr(A_1)=0$, and zeroes elsewhere; all $A_i=0$ gives $N^+$ back. 
Explicitly, then, a general member of the slice has the shape
\[
A = \begin{pmatrix} A_1 & I \\ A_2 & & I \\ \vdots \\ A_{m-1} &&&&& I
\\ A_m && \cdots &&& 0 \end{pmatrix}
\]
where the $A_k$ are $2 \times 2$ matrices, and with $tr(A_1) = 0$.
The characteristic polynomial is $det(\lambda I - A) = det( \lambda^m
- A_1 \lambda^{m-1} - \cdots - A_m)$. The smooth fibres $Y_{m,t} =
\chi^{-1}(t) \cap \Slice_m$ are smooth complex affine varieties of
dimension $2m$.

Caution: for Lie theory purists, this is not in fact a JM slice (there is no suitable $N^-$), but is orbit-preservingly isomorphic to JM slices, and more technically convenient for our purposes.

Note that the generic fibre of the map (i.e. over a point of configuration space) is unchanged by the simultaneous resolution, so in principle all the topology of these spaces can be understood explicitly in terms of the linear algebra of certain matrices.  On the other hand, the resolution of the zero-fibre (the nilpotent cone) is well-known to retract to a compact core which is just the preimage of the matrix $N^+$ itself; in other words, it's the locus of all flags stabilised by $N^+$. This core is a union of \emph{complex Lagrangian} submanifolds, described in more detail in \cite{Kho:CM}, in particular the number of irreducible components is given by the Catalan number $\frac{1}{m+1}{2m \choose m}$.

In nearby smooth fibres, these complex Lagrangian components $L_{\wp}$ give rise to distinguished real Lagrangian submanifolds, and it is plausible to conjecture that this finite set of Lagrangian submanifolds generate Donaldson's quantum category of $Y_{m,t}$ (the underlying homological category of the Fukaya category) in the weak sense that every Floer homologically essential closed Lagrangian submanifold has non-trivial Floer homology with one of the $L_{\wp}$, cf. Section 1(J).
\medskip

\noindent \textbf{2(g) Inductive geometry:}
The key construction with this slice is an ``inductive scheme'', relating the ``least singular non-smooth'' fibres of $\Slice_m$ to the smooth fibres of $\Slice_{m-1}$.  Fix $\mu = (\mu_1=0,\mu_2=0, \mu_3,\ldots,\mu_m)$ a tuple of eigenvalues, the first two of which vanish and with all others being pairwise distinct, and let $\hat{\mu} = (\mu_3,\ldots,\mu_m)$.  

Lemma: The fibre of $(\chi|\Slice_m)^{-1}(\mu)$ has complex codimension 2 smooth singular locus which is canonically isomorphic to $(\chi|\Slice_{m-1})^{-1}(\hat{\mu})$. Moreover, along the singular locus $\chi$ has a fibred $A_1$-singularity (an open neighbourhood of the singular locus looks like its product with $x^2+yz=0$).

Rescaled parallel transport and the vanishing cycle construction give a relative vanishing cycle in smooth fibres of $\chi|\Slice_m$ which is an $\S^2$-bundle over a fibre of $\chi|\Slice_{m-1}$.  General properties of symplectic parallel transport give that these relative vanishing cycles are not Lagrangian but coisotropic, with the obvious $\S^2$-fibrations being the canonical foliations by isotropic leaves.

The force of the Lemma is that this process can now be \emph{iterated}.  Of course, an isotropic fibration restricted to a Lagrangian submanifold gives rise to a Lagrangian submanifold of the total space.
\medskip

\noindent \textbf{2(h) Symplectic Khovanov homology:}
Fix a crossingless matching $\wp$ of $2m$ points in the plane; the points specify a fibre $Y_m$ of $\chi|\Slice_m$. Bringing eigenvalues together in pairs along the paths specified by the matching, and iterating the vanishing cycle construction above, gives a Lagrangian $L_{\wp}$ which is an iterated $\S^2$-bundle inside $Y_m$.  In fact one can show that it is diffeomorphic to $(\S^2)^m$ (hence spin). We care especially about the first case $\wp_+$ below; we remark that the number of crossingless matchings which lie entirely in the upper half-plane, up to isotopy, is given by the Catalan number $\frac{1}{m+1}{2m \choose m}$, cf. section (f) above.

Crossingless matchings: the nested horseshoe on the left is denoted $\wp_+$.

\[
 \includegraphics{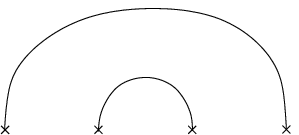}
 \qquad \qquad
 \includegraphics{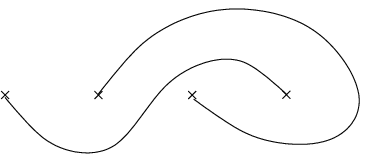}
\]

Given a braid $\beta \in Br_m$ we get a Floer group via thinking of $\beta \times \id \in Br_{2m} = \pi_1(\Conf_{2m}(\C))$ as explained above.  The following is the main result of these notes:

\begin{Theorem}[Seidel, S.] \label{Theorem} The $\Z$-graded Floer cohomology group 
$$Kh_{symp}^*(K_{\beta}) = HF^{*+m+w}(L_{\wp_+},(\beta\times\id)(L_{\wp_+}))$$ 
is an oriented link invariant: here $m$ is the number of strands and $w$ the writhe of the braid diagram.
\end{Theorem}

It is important to realise that the loss of information in passing from the bigrading to the single grading is substantial: for instance, $Kh_{symp}$ does not in itself determine the Jones polynomial. 

The proof of the Theorem involves verifying invariance of the Floer group under the Markov moves. For the first move, this is relatively straightforward, since -- once the machinery of rescaled parallel transport has been carefully set in place -- the Lagrangian submanifold $L_{\wp}$ is itself unchanged (up to Hamiltonian isotopy) by effecting a conjugation.  For the second Markov move, the proof is more involved since one must compare Floer groups for Lagrangians of different dimension living in different spaces.  The key is the fibred $A_1$-structure along the singular set where two eigenvalues co-incide, and a fibred $A_2$-generalisation to the case where three eigenvalues co-incide. Indeed, locally near branches of the discriminant locus of $\chi$ where eigenvalues $1,2$ resp. $2,3$ co-incide, the smooth fibres contain pairs of vanishing cycles which together form an $\S^2\vee\S^2$-fibration.  The fact that the monodromy symplectomorphisms about the two branches of the discriminant satisfy the braid relations can be deduced explicitly from Section 1(H); of course, the fact that we have a fibre bundle over configuration space gives the same result without any appeal to the local structure.

However, a similar local analysis allows one to explicitly identify the Floer complexes for the Lagrangians before and after the Markov $II$ move.  (The grading shift in the definition takes care of the difference of the effects of the Markov $II^+$ and Markov $II^-$ moves on the Maslov class.)  The upshot is that very general features of the singularities of the mapping $\chi$ encode the local geometric properties which lead to the symplectic Khovanov homology being an invariant. 
\medskip

\noindent \textbf{2(i) Long exact sequences revisited:}
In a few cases -- unlinks, the trefoil -- one can compute $Kh_{symp}$ explicitly, and in such cases one finds that the answer agrees with Khovanov's combinatorial theory.  Even in these cases, the result is rather surprising, since the methods of computation do not particularly parallel one another.  Thus the Main Theorem is complemented by:

\begin{Conjecture}[Seidel, S.] \label{Conjecture} $Kh^*_{symp} = \oplus_{i-j=*}Kh^{i,j}$.
\end{Conjecture}

Main evidence: $Kh^*_{symp}$ should also satisfy the right skein-type exact sequences (in the notation of Section 2(b))
\begin{small}
$$Kh^*_{symp}(L_{hor})\rightarrow Kh^*_{symp}(L_{vert})\rightarrow Kh^*_{symp}(L_{cross}).$$
\end{small}
\noindent These should come from a version of the LES in Floer theory for a \emph{fibred} Dehn twist, which is just the monodromy of $\chi$ corresponding to inserting a single negative crossing.  Indeed, one can speculate that appropriate long exact sequences exist for suitable correspondences, as follows.  

Suppose in general we are given Lagrangians $L_0,L_1 \subset X$ and $\hat{L}_0, \hat{L}_1 \subset \hat{X}$, and a Lagrangian correspondence $C \subset (X\times \hat{X}, \omega_X\oplus -\omega_{\hat{X}})$ which is an isotropic $\S^a$-fibration over $X$.  Suppose moreover the $\hat{L}_i$ are given by lifting the $L_i$ from $X$ to $\hat{X}$ via the correspondence.  One can try to find an exact triangle of the shape
$$HF(C\times C, L_0 \times \Delta \times L_1) \rightarrow HF(\hat{L}_0, \hat{L}_1) \rightarrow HF(\hat{L}_0, \tau(\hat{L}_1))$$
where $\tau$ denotes a fibred Dehn twist along $C$ and the first homology group is taken inside $X\times\hat{X}\times\hat{X}\times{X}$ with the symplectic form reversed on the second two factors, and with $\Delta$ the diagonal.  Moreover, if the geometry is sufficiently constrained, one can hope to relate the first group to Floer homology $HF(L_0,L_1)$ taken inside $X$.  Using the relative vanishing cycles inside $\Slice_{m-1}\times\Slice_m$ as correspondences, this general picture includes the desired skein-type relation.

If one assumes the existence of the long exact sequence, then following a rather general algebraic strategy one can construct a spectral sequence with $E^2=Kh^{*,*}$ and converging to $E^{\infty}=Kh^*_{symp}$ (the model outline is contained in Ozsvath and Szabo's work \cite{OzSz:DBC}, in which they use a similar approach to relate Khovanov homology of a link $L$ with the Heegaard Floer homology of the branched double cover $M(L)$).  From this perspective, the above conjecture asserts the vanishing of the higher differentials in this spectral sequence; in the analogous story with Heegaard Floer theory, by contrast, the higher order differentials are often non-zero.

A distinct circle of ideas relating the chain complex underlying symplectic Khovanov homology to the Bigelow-Lawrence homological construction of the Jones polynomial \cite{BigLaw} has recently been given by Manolescu in \cite{Man}.
\medskip

\noindent \textbf{2(j) Counting sections revisited:}
Khovanov's theory is especially interesting since it fits into a TQFT (and we know a lot about knots, but little about surfaces with or without boundary in $\R^4$).  A small piece of that is easily visible in $Kh^*_{symp}$, in the spirit of Section~1(G);~(2). 

Suppose we have a \emph{symplectic} cobordism (surface in $\R^4$) between two positive braids. By fibring $\R^4 \subset \C\P^2$ by $\C$-lines, we get a braid monodromy picture of the surface, which is just a relative version of the Lefschetz fibration story from Section 1(F).  Geometrically, the braid monodromy gives  an annulus in configuration space whose boundary circles represent the two boundary knots/braids.  Now counting holomorphic sections of $\chi$ over the annulus, with suitable Lagrangian boundary conditions, gives rise to a morphism on Floer homology groups and hence on symplectic Khovanov homology.

Challenge: detect symplectically knotted surfaces (or families of such with common boundary) this way.

\end{document}